\documentclass[12pt]{article}
\usepackage{amssymb}

\newcommand{\F}{\mathcal{F}}
\newcommand{\G}{\mathcal{G}}
\newcommand{\A}{\mathcal{A}}
\newcommand{\B}{\mathcal{B}}
\newcommand{\Oh}{\mathcal{O}}
\newcommand{\La}{{\rm La}}
\newcommand{\E}{{\rm E}}
\renewcommand{\Pr}{{\rm Pr}}
\newtheorem{lemma}{Lemma}
\newtheorem{theorem}{Theorem}

\def\nchn{{n \choose \lfloor \frac{n}{2}\rfloor}}
\def\uppernk{{\lfloor \frac{n+(k-2)}{2} \rfloor}}
\def\lowernk{{\lfloor \frac{n-(k-2)}{2} \rfloor}}
\def\lanh{\La(n,H)}
\newcommand{\comments}[1]{}
\begin{document}
\title{On families of subsets with a forbidden subposet}
\author{Jerrold R. Griggs\thanks{
Department of Mathematics, University of South Carolina, Columbia,
SC 29208  USA ({\tt griggs@math.sc.edu}).  This author's work was
supported in part by the conference ``New Directions in
Algorithms, Combinatorics, and Optimization," held in honor of Tom
Trotter at Georgia Tech in May, 2008.} \and Linyuan Lu
\thanks{Department of Mathematics,
University of South Carolina, Columbia, SC 29208  USA ({\tt
lu@math.sc.edu}). This author was supported in part by NSF grant
DMS 0701111. } }

\date{ June 7, 2008\\}

\maketitle

\begin{abstract}
Let $\F\subset 2^{[n]}$ be a family of subsets of $\{1,2,\ldots,
n\}$. For any poset $H$, we say $\F$ is $H$-free if $\F$ does not
contain any subposet isomorphic to $H$. Katona and others have
investigated the behavior of $\La(n,H)$, which denotes the maximum
size of $H$-free families $\F\subset 2^{[n]}$. Here we use a new
approach, which is to apply methods from extremal graph theory and
probability theory to identify new classes of posets $H$, for
which $\La(n,H)$ can be determined asymptotically as
$n\rightarrow\infty$ for various posets $H$, including
two-end-forks, up-down trees, and cycles $C_{4k}$ on two levels.
\end{abstract}

\emph{Dedicated to Prof. William T. Trotter on the occasion of his
65th birthday}

\section{Introduction and Results}
A poset $(S,\leq)$ is a set $S$ equipped with a partial ordering
$\leq$. We say a poset $(S,\leq)$ contains another poset
$(S',\leq')$ if there exists an injection $f\colon S' \to S$,
which preserves the partial ordering, meaning that whenever
$u,v\in S'$ satisfy $u\leq' v$, we have $f(u)\leq f(v)$. In this
case, $S'$ is called a subposet of $S$.

Let $\F\subset 2^{[n]}$ be a family of subsets of
$[n]:=\{1,2,\ldots, n\}$. For any poset $H$, we say $\F$ is
$H$-free if the poset $(\F,\subseteq)$ does not contain $H$ as a
subposet. Let $\La(n,H)$ denote the largest size of $H$-free
family of subsets of $[n]$.  The fundamental result of this kind
is for $H$ being a chain $P_2$ of two elements.   A $P_2$-free
family is an antichain, and Sperner's Theorem~\cite{sperner} from
1928 gives us that $\La(n,P_2)=\nchn$. For small posets $H$ in
general, it is interesting to compare $\La(n,H)$ to $\nchn$.

\comments{Here we use a new approach, which is to apply methods
from extremal graph theory and probability theory to identify new
classes of posets $H$, for which $\La(n,H)$ can be determined
asymptotically.

 Sperner's Theorem can be derived from
Lubell-Yamamoto-Meshalkin's inequality
\begin{equation}
  \label{eq:LYM}\[  \sum_{F\in \F}\frac{1}{{n\choose |F|}}\leq 1 \]
\end{equation}
for any antichain $\F$. }

Erd\H os~\cite{erdos} extended Sperner's Theorem in 1945 to
determine that $\La(n,P_k)$, where $P_k$ is a chain (path) of $k$
elements, is the sum of the $k-1$ middle binomial coefficients in
$n$.  Consequently, $\La(n,P_k)\sim(k-1)\nchn$, as
$n\rightarrow\infty$.  Let $h(P)$ denote the height of poset $P$,
which is the largest cardinality of any chain in $H$.  We are
interested in the asymptotic behavior of $\lanh$ for other posets
$H$ of {\it height\/} $k$.

There have been several investigations already of height two
posets.  Thanh~\cite{thanh} extended Sperner's Theorem by showing
that for all $r$, $\La(n,V_r)\sim\nchn$, where $V_r$ is the
$r$-fork, the height two poset with one element at the bottom
level below each of $r$ elements at the top level.  (Especially,
$V_1$ is $P_2$, while $V_2$ looks like the letter $V$.)  It is
important to note that we are not only excluding ``induced" copies
of a forbidden subposet $H$, e.g., $V_3$ is a subposet of $P_4$,
so excluding $V_3$ subposets also excludes $P_4$.

DeBonis and Katona~\cite{dk} determined that $\La(n,B)$, where $B$
is the Butterfly poset on four elements $A_1,A_2,B_1,B_2$ with
each $A_1,A_2\le B_1,B_2$, is the sum of the two middle binomial
coefficients in $n$.  More generally, consider excluding the
height two poset which is called (using graph-theoretic
terminology) $K_{r,s}$, which has elements $A_i, 1\le i\le r$ at
the bottom level, elements $B_j, 1\le j\le s$ at the top level,
and for all $i,j$, $A_i\le B_j$.  DeBonis and Katona\cite{dk}
extend the asymptotics for the butterfly $B$ and show that
$\La(n,K_{r,s})\sim2\nchn$ for all $r,s\ge2$.  Griggs and
Katona~\cite{gk} considered whether the asymptotics of excluding
the $N$ poset on four elements $A_1,A_2,B_1,B_2$ with $A_1\le B_1,
A_2\le B_1, A_2\le B_2$ is similar to excluding $V_2$ or $B$. It
turns out to be the former: $\La(n,N)\sim \nchn$.

One new class of posets considered here we call a {\it baton\/}
$P_k(s,t)$, which is a path $P_k$ on $k$ elements, $k\ge3$, such
that the bottom element is replicated $s-1$ times and the top
element is replicated $t-1$ times, $s,t\ge1$.  That is, we have a
height $k$ poset with $s$ (resp. $t$) independent elements on the
bottom (resp., top) level. The particular case $P_k(1,r)$ (which
resembles a palm tree), known as an $r$-fork with a $k$-shaft, has
been examined by Katona and De Bonis \cite{dk}. They show
\begin{eqnarray}
\label{eq:2} La(n,P_k(1,r))&\geq& \sum_{i=\lowernk} ^{\uppernk}
{n\choose i} + {n\choose \lfloor \frac{n+k+1}{2}\rfloor}
 \left(\frac{r-1}{n}
+ \Omega(\frac{1}{n^2})\right)\\
\label{eq:3} La(n,P_k(1,r))&\leq&
 \sum_{i=\lowernk}
^{\uppernk} {n\choose i} + {n\choose \lfloor
\frac{n+k+1}{2}\rfloor}
 \left(\frac{z(k)+2(r-1)}{n}
+ \Omega(\frac{1}{n^2})\right)
\end{eqnarray}
where $z(k)=\lfloor \frac{k^2}{2}\rfloor$ if $n+k$ is even
and $z(k)=\lfloor \frac{(k-1)^2}{2}\rfloor$ if $n+k$ is odd.

The previously known maximum sizes of families of subsets of $[n]$
without a given pattern are listed in the following table.
\vfill\eject

\begin{table}[htbp]
  \centering
  \begin{tabular}{|c|c|c|c|}
  \hline
\hline
Name & $H$ & $\La(n, H)$ & Reference\\
\hline Chain $P_r$ & $A_1\subset \cdots \subset A_r$
&$(r-1+o_n(1)){[n] \choose \lfloor \frac{n}{2}\rfloor}$
& \cite{erdos}\\
\hline Butterfly $B$ &$A_i \subset B_j$, for $1\leq i,j \leq 2$
&$(2+o_n(1)){[n] \choose \lfloor \frac{n}{2}\rfloor}$
 & \cite{dks}\\
\hline
$K_{r,s}$ ($r,s\geq 2$) &$A_i \subset B_j$, for $1\leq i\leq r$,
$1\leq j \leq s$
&$(2+o_n(1)){[n] \choose \lfloor \frac{n}{2}\rfloor}$
 & \cite{dk}\\
\hline
``N'' & $A\subset B$,  $C\subset B$, and $C\subset D$
&$(1+o_n(1)){[n] \choose \lfloor \frac{n}{2}\rfloor}$
 & \cite{gk}\\
\hline
``$V_r$'' & $A\subset B_i$, for $i=1,2,\ldots, r$ &
$(1+o_n(1)){[n] \choose \lfloor \frac{n}{2}\rfloor}$ & \cite{thanh}\\
\hline
$_kV_r$ & $A_1\subset\cdots\subset A_k \subset B_i$,
for $i=1,2,\ldots, r$ &
$(k+o_n(1)){[n] \choose \lfloor \frac{n}{2}\rfloor}$ & \cite{dk}\\
\hline
  \end{tabular}
  \caption{Previously known results in the literature}
  \label{tab:1}
\end{table}

In this paper we give new asymptotic upper bounds on
$\La(n,H)/\nchn$ for several classes of posets $H$, and identify
some new ones for which this ratio goes to 1 as
$n\rightarrow\infty$.  We first "roughly unify" the previous
results on forks $_kV_r$ and on complete two level posets
$K_{s,t}$ by considering batons $P_k(s,t)$. Note that the
summation term in the bound, which appears repeatedly, is just the
sum of the $k-1$ middle binomial coefficients in $n$.

\begin{theorem}\label{t1}
  For any $s,t\geq 1$ and $k\geq 3$,
We have
\begin{equation}
  \label{eq:1}
  \La(n,P_k(s,t))\leq
 \sum_{i=\lowernk}^{\uppernk} {n\choose i} +
 {n\choose \lfloor\frac{n+k}{2}\rfloor}
 \left(\frac{2k(s+t-2)}{n}
+ O(n^{-3/2}\sqrt{\ln n})\right).
\end{equation}
Consequently, as $n\rightarrow\infty$,
$$
\La(n,P_k(s,t))/\nchn\rightarrow k-1.
$$
\end{theorem}

\noindent {\bf Remarks:}
\begin{enumerate}

\item
Theorem \ref{t1} (at $s=1$ and $t=r$) is better than inequality
(\ref{eq:3}) for $k\geq 4r-3$. For small $k$ and large $r$,
inequality (\ref{eq:3}) gives a better constant in the second
order term.
\item Note $\La(n,P_k(s,t))\geq \La(n,P_k(1,\max\{s,t\}))$.
>From inequality (\ref{eq:2}), we have
\begin{equation}
  \label{eq:4}
 La(n,P_k(s,t))\geq
\sum_{i=\lowernk} ^{\uppernk} {n\choose i} + {n\choose \lfloor
\frac{n+k}{2}\rfloor}
 \left(\frac{\max\{s,t\}-1}{n}
+ \Omega(\frac{1}{n^2})\right).
\end{equation}
This lower bound (\ref{eq:4}) can be compared to the upper bound
(\ref{eq:1}).
\item Note that $P_3(s,t)$ contains $P_2(s,t)=K_{s,t}$,
the complete two level poset.
Theorem \ref{t1} implies
\begin{equation}
  \label{eq:5}
\La(n,H)\leq \left(2+O(\frac{|H|}{n})\right) {\nchn}
\end{equation}
for all posets of height $2$. The hidden constant in the second
order term is slightly worse than that given in \cite{dk}. If $H$
is not a subposet of the two middle layers of $2^{[n]}$ (for
example $H$ contains the butterfly $B$), then the equality in
(\ref{eq:5}) holds.
\end{enumerate}

An {\em up-down tree\/} $T$ is a poset of height $2$ that is also
a tree as an undirected graph; its order is the number of
elements, $|T|$.
\begin{theorem} \label{t3}
  For any up-down tree $T$ with order $t$, we have
  \begin{equation}
    \label{eq:6}
  \La(n, T)\leq \left(1+ \frac{16t}{n} +
  O\left(\frac{1}{n\sqrt{n\ln n}}\right)\right){\nchn}.
  \end{equation}
Consequently, as $n\rightarrow\infty$,
$$
\La(n,T)/\nchn\rightarrow 1.
$$
\end{theorem}

After discovering the results above for batons and for up-down
trees, we learned of new progress by Boris Bukh~\cite{bukh} that
describes the asymptotic behavior of $\La(n,T)$ for every tree
poset. Specifically, if $T$ is any poset for which the Hasse
diagram is a tree (connected and acyclic), then
\begin{equation}
\La(n,T)=(h(T)-1)\nchn (1+O(1/n)).
\end{equation}
This implies the leading asymptotic behavior for batons and
up-down trees in Theorems \ref{t1} and \ref{t3} above, though the
proofs and error terms are different.

The butterfly poset $B$ has been solved asymptotically, so it is
next interesting to consider more generally the crowns $\Oh_{2k}$,
which is the poset of height $2$ that is a cycle of length $2k$ as
an undirected graph. Of course, $\Oh_4$ is the butterfly poset,
while $\Oh_6$ is noteworthy for being the middle two levels of the
Boolean lattice $B_3$.  We have the following theorem for crowns:

\begin{theorem} \label{t4}
  For $k\geq 2$, we have
  \begin{eqnarray}
\label{eq:7} \La(n, \Oh_{4k})&=& (1+o_n(1))
{\nchn}\\
\label{eq:8} \La(n, \Oh_{4k-2})&\leq& \left(1+\frac{\sqrt{2}}{2}+
o_n(1)\right) {\nchn}.
  \end{eqnarray}
\end{theorem}

So we see that the crowns $\Oh_{2k}$, $k\ge3$, have
$\La(n,\Oh_{2k})/\nchn$ staying strictly below 2 asymptotically,
unlike the Butterfly, the case $k=2$, where the ratio goes to 2.
For even $k\ge4$, the ratio goes to 1, while for odd $k\ge3$ we
only have an asymptotic upper bound.

The Theorem above for crowns is actually just a special case of
the more general result which concerns a more general class of
height 2 posets obtained from graphs in a natural way.  The proof
also relies on extremal graph theory.  For a simple graph
$G=(V,E)$, define a poset $P(G)$ on the set $V \cup E$ with the
partial ordering $v<e$ if the edge $e$ is incident at vertex $v$
in $G$.  For example, the crown poset $\Oh_{2k}$ is $P(G)$ when
graph $G$ is a $k$-cycle.

\begin{theorem} \label{t5}
For any nonempty simple graph $G$ with chromatic number $\chi(G)$,
we have
\begin{equation}
  \label{eq:9}
\La(n, P(G))\leq \left(1+ \sqrt{1-\frac{1}{\chi(G)-1}}+
o_n(1)\right) {\nchn}.
\end{equation}
In particular, if $G$ is a bipartite graph, then
\begin{equation}
  \label{eq:10}
\La(n, P(G))= (1+o_n(1)) {\nchn}.
\end{equation}
\end{theorem}

Theorem \ref{t4} is a direct consequence of Theorem \ref{t5} by
the observation $\Oh_{2k}=P(C_{k})$.

In this theory we construct large families in the Boolean lattice
that avoid a given subposet.  This is analogous to the
much-studied Tur\'an theory of graphs, in which one seeks to
maximize the number of edges on $n$ vertices while avoiding a
given subgraph. It is interesting that the theorem above applies
the Tura\'an theory of graphs to give a useful bound in our
ordered set theory.

The rest of the paper is organized as follows. Three probabilistic
lemmas are given in Section 2, and the proofs of the theorems are
given in section 3.  We conclude with ideas for further research.

\section{Lemmas}
For any fixed poset $H$, $\La(n,H)$ is of  magnitude
$\Theta\left({\nchn}\right)$. The following lemma allows us to
consider the families consisting only of subsets near the middle
level.

\begin{lemma}\label{l1}
For any positive integer $n$, we have
  \begin{eqnarray}
\label{eq:11}
  \sum_{i>\frac{n}{2}+ 2\sqrt{n\ln n}} {n\choose i}&<& \frac{2^n}{n^2};\\
\label{eq:12}
  \sum_{i<\frac{n}{2}- 2\sqrt{n\ln n}} {n\choose i}&<& \frac{2^n}{n^2}.
  \end{eqnarray}
\end{lemma}
{\bf Proof}: Let $X_1, X_2,\ldots, X_n$ be $n$ independent identically
distributed $\{0,1\}$ random variables with
$$\Pr(X_i=0)=\Pr(X_i=1)=\frac{1}{2}$$
for any $1\leq i\leq n$. Apply Chernoff's inequality
\cite{chernoff} to $X=\sum_{i=1}^n X_i$. We have
\[
  \Pr(X-E(X)>\lambda)< e^{-\frac{\lambda^2}{2n}}.
\]
Choose $\lambda=2\sqrt{n\ln n}$. We have
\begin{eqnarray*}
  \sum_{i>\frac{n}{2}+ 2\sqrt{n\ln n}} {n\choose i}2^{-n}
&=& \Pr(X>\frac{n}{2}+\lambda)\\
&<&e^{-\frac{\lambda^2}{2n}}\\
&=&\frac{1}{n^2}.
\end{eqnarray*}
Inequality (\ref{eq:11}) has been proved.
Inequality (\ref{eq:12}) is equivalent to inequality (\ref{eq:11})
by the symmetry of binomial coefficients
${n\choose i}={n\choose n-i}$.
\hfill $\square$

Apply Stirling's formula $n!=(1+O(1/n)\sqrt{2\pi n})
\frac{n^n}{e^n}$ to obtain the following approximation of
${n\choose \lfloor \frac{n}{2}\rfloor}$:
\begin{eqnarray*}
  {n\choose \lfloor \frac{n}{2}\rfloor}&=&
 \frac{n!}{\lfloor \frac{n}{2}\rfloor ! \lceil \frac{n}{2}\rceil !}\\
&=&(1+O(1/n)) \frac{\sqrt{2\pi n}{\frac{n^n}{e^n}}}
{\sqrt{2\pi\lfloor \frac{n}{2}\rfloor }{\frac{\lfloor
\frac{n}{2}\rfloor^{\lfloor \frac{n}{2}\rfloor}}{e^{\lfloor
\frac{n}{2}\rfloor}}} \sqrt{2\pi \lceil
\frac{n}{2}\rceil}{\frac{\lceil \frac{n}{2}\rceil
^{\lceil \frac{n}{2}\rceil}}{e^{\lceil \frac{n}{2}\rceil}}}}\\
&=& (1+O(1/n))\frac{\sqrt{2}}{\sqrt{\pi n}} 2^n.
\end{eqnarray*}

It implies that $\frac{2^n}{n^2}=(1+O(1/n))
\frac{\sqrt{\pi/2}}{n^{3/2}} {n\choose \lfloor
\frac{n}{2}\rfloor}$. For any family $\F$ of size
$\Theta({n\choose \lfloor \frac{n}{2}\rfloor})$, we can delete all
subsets of sizes not in $(\frac{n}{2}- 2\sqrt{n\ln n},
\frac{n}{2}+ 2\sqrt{n\ln n})$ from $\F$. We obtain a family of
subsets that has about the same size of $\F$ and only contains
subsets of sizes in $(\frac{n}{2}- 2\sqrt{n\ln n}, \frac{n}{2}+
2\sqrt{n\ln n})$.



\begin{lemma}\label{l4}
Suppose $X$ is a random variable which takes on nonnegative
integer values. Let $f(x)$ and $g(x)$ be two nondecreasing
functions defined for nonnegative integers $x$. Then
\[
\E(f(X)g(X))\geq \E(f(X)) \E(g(X)).
\]
\end{lemma}
{\bf Proof:} Apply the FKG inequality~\cite{fkg} over the totally
ordered set of nonnegative integers.  Alternately, here we give a
simple direct proof.

For any integer $k\geq 1$, let $h_k$ be the step function:
\[
h_k(x)=\left\{
  \begin{array}{ll}
0 & \mbox{ if } 0\leq x <k;\\
1 & \mbox{ if } x \geq k .
  \end{array}
\right.
\]
For integers $j\ge i\ge1$, we observe that
\[\E(h_i(X)h_j(X))\geq \E(h_i(X))\E(h_j(X)),\]
which holds sine
\begin{eqnarray*}
  \E(h_i(X)h_j(X))&=& \Pr(X\geq i\ \&\ X\geq j)\\
 &=& \Pr(X \geq j)\\
 &\geq&\Pr(X \geq i) \Pr(X \geq j)\\
 &=&\E(h_i(X))\E(h_j(X)).
\end{eqnarray*}
We have
\[f(x)= f(0) + \sum_{k=1}^\infty (f(k)-f(k-1)) h_k(x).\]
Similarly
\[g(x)=g(0) +  \sum_{k=1}^\infty (g(k)-g(k-1)) h_k(x).\]
All coefficients $f(k)-f(k-1)$ and $g(k)-g(k-1)$ are nonnegative.
By linearity, we have
\begin{eqnarray*}
\E(f(X)g(X)) &=& \E( (f(0) + \sum_{i=1}^\infty (f(i)-f(i-1))
h_i(X))
(g(0) +  \sum_{j=1}^\infty (g(j)-g(j-1)) h_j(X)))\\
&=& f(0)g(0) + f(0) \sum_{j=1}^\infty (g(j)-g(j-1))\E(h_j(X))\\
&&+ g(0)  \sum_{i=1}^\infty (f(i)-f(i-1)) \E(h_i(X))\\
&&+ \sum_{i,j=1}^\infty (f(i)-f(i-1))(g(j)-g(j-1))\E(h_i(X)h_j(X))\\
&\geq&
 f(0)g(0) + f(0) \sum_{j=1}^\infty (g(j)-g(j-1))\E(h_j(X))\\
&&+ g(0)  \sum_{i=1}^\infty (f(i)-f(i-1)) \E(h_i(X))\\
&&+ \sum_{i,j=1}^\infty (f(i)-f(i-1))(g(j)-g(j-1))\E(h_i(X))\E(h_j(X))\\
&=& \E(f(X))\E(g(X)).
\end{eqnarray*}
 \hfill $\square$

\begin{lemma}\label{l3}
Suppose $X$ is a random variable which takes on nonnegative
integer values.  For integers $k>r\geq 1$, if $\E(X)> k-1$, then
\begin{equation}
  \label{eq:13}
\E{X\choose k}\geq \E{X\choose r}
\frac{r!}{k!}\prod_{i=0}^{k-r-1}(\E(X)-r-i).
\end{equation}
\end{lemma}
{\bf Proof:}
 Define \[f(x)=\left\{
    \begin{array}[c]{ll}
\frac{r!}{k!}\prod_{i=0}^{k-r}(x-r-i)
& \mbox{ if }x>k-1\\
0 &\mbox{ otherwise.}
    \end{array}
\right.
\]
and
\[g(x)=\left\{
    \begin{array}[c]{ll}
\frac{1}{r!}\prod_{i=0}^{r-1}(x-i)
& \mbox{ if } x > r-1\\
0 &\mbox{ otherwise.}
    \end{array}
\right. .\] Both $f(x)$ and $g(x)$ are nonnegative increasing
functions.
For each nonnegative integer $x$, we have $g(x)={x\choose r}$ and
$f(x)g(x)={x\choose k}$. By applying Lemma~\ref{l4} we obtain
\begin{eqnarray*}
  \E{X\choose k} &=& \E(f(X)g(X))\\
  &\geq& \E(f(X)) \E(g(X))\\
  &=& \E(f(X)) \E{X\choose r}\\
  &\geq& f(\E(X))\E{X\choose r},
\end{eqnarray*}
where the last inequality follows from since $f(x)$ is concave
upward. \hfill$\square$
\section{Proofs of theorems}
\noindent {\bf Proof of Theorem \ref{t1}:} We let
$\epsilon=\frac{2k(s+t-2)}{\frac{n}{2}-2\sqrt{n\ln n}},$ and
\[f=f(n,k,s,t)=\sum_{i=\lowernk}
^{\uppernk} {n\choose i} + {n\choose \lfloor
\frac{n+k}{2}\rfloor}\epsilon.\] Suppose $\F$ is a family of
subsets of $[n]$ with $|\F|>f +\frac{2^{n+1}}{n^2}$. By removing
all subsets of size outside $(\frac{n}{2}- 2\sqrt{n\ln n},
\frac{n}{2}+ 2\sqrt{n\ln n})$, we can assume $\F$ only contains
subsets of sizes in $(\frac{n}{2}- 2\sqrt{n\ln n}, \frac{n}{2}+
2\sqrt{n\ln n})$ and $|\F|> f$.

We would like to show $\F$ contains $P_k(s,t)$. We will prove this
statement by contradiction. Suppose that $\F$ is $P_k(s,t)$-free.
Take a random permutation $\sigma\in S_n$. Consider a random full
(maximal) chain $C_\sigma$
\[\emptyset\subset \{\sigma_1\} \subset  \{\sigma_1,\sigma_2\}
\subset \cdots \subset \{\sigma_1,\sigma_2,\cdots, \sigma_n\}.\]
Let $X$ be the random number counting $|\F \cap C_\sigma|$. On the
one hand, we have
\begin{eqnarray}
\label{eq:15}
  \E(X)&=& \sum_{F\in \F}\frac{1}{{n\choose |F|}}\\
\label{eq:16} &>& k-1 +\epsilon,
\end{eqnarray}
since the sum is minimized, for a family of subsets on $[n]$ of
size $f$ by taking the $f$ sets closest to the middle size $n/2$,
which means taking the $k-1$ middle levels and the remaining sets
at the next closest level to the middle, $\lfloor
\frac{n+k}{2}\rfloor$.

Apply Lemma \ref{l3} with $r=k-1$
\begin{eqnarray}
\nonumber
\E{X\choose k}&\geq&
\frac{1}{k}\E{X\choose k-1}(\E(X)-k+1)\\
\label{eq:17}
&>&\frac{\epsilon}{k}\E{X\choose k-1}.
\end{eqnarray}

On the other hand, we will compute $\E{X\choose k}$ directly. By
coumting chains, a subchain of length $k$ in $\F$,
\[F_1\subset F_2\subset \cdots \subset F_{k},\]
is in the random chain $C_\sigma$ with probability
\[\frac{|F_1|! (|F_2|-|F_1|)!\cdots (n-|F_{k}|)!}{n!}.\]
By linearity, we have
\begin{equation}
\label{eq:18}
\E{X\choose k}=\sum_{\stackrel{F_1,\ldots, F_{k}\in \F}
{F_1\subset \cdots \subset F_{k}}}
\frac{|F_1|! (|F_2|-|F_1|)!\cdots (n-|F_{k}|)!}{n!}.
\end{equation}

We can rewrite equation (\ref{eq:18}) as
\begin{equation}
  \label{eq:19}
\E{X\choose k}= \!\!\!\!\!\!
\sum_{\stackrel{F_{2}, \ldots, F_{k-1}\in \F}
{F_{2}\subset \cdots \subset F_{k-1}}}
\!\!\!\!\!\!\!
\frac{|F_{2}|! \cdots (n-|F_{k-1}|)!}{n!}
\sum_{\stackrel{F_1\in \F}{F_1\subset F_2}} \frac{1}{{|F_2| \choose |F_1|}}
\!\!\!\!
\sum_{\stackrel{F_{k}\in \F}{F_{k-1}\subset F_{k}}}
\frac{1}{{n-|F_{k-1}| \choose n-|F_{k}|}}.
\end{equation}

Since $\F$ is $P_k(s,t)$-free, for a fixed
$F_2,\ldots, F_{k-1}$,
either ``the number of $F_1$ satisfying $F_1\subset F_2$ is at most
 $s-1$''
or ``the number of $F_{k}$ satisfying $F_{k-1}\subset F_{k}$ is at
most $t-1$''. Let $\A$ be the set of $k-2$-chains $F_2\subset
\ldots \subset F_{k-1}$ in $\F$ so that the number of $F_1\in \F$,
$F_1\subset F_2$, is at most
 $s-1$. Let $\B$ be the set of $k-2$-chains $F_2\subset \ldots \subset F_{k-1}$
 in $\F$ so that
the number of $F_k\in \F$, $F_{k-1}\subset F_k$, is at most
 $t-1$. The union of $\A$ and $\B$ covers all $k-2$-chains in $\F$.
We have
\begin{eqnarray}
\nonumber
  \E{X\choose k} &\leq& \!\!\!\!\!\!
\sum_{(F_{2}, \ldots, F_{k-1})\in \A}
\!\!\!\!\!\!\!
\frac{|F_{2}|! \cdots (n-|F_{k-1}|)!}{n!}
\sum_{\stackrel{F_1\in \F}{F_1\subset F_2}} \frac{1}{{|F_2| \choose |F_1|}}
\!\!\!\!
\sum_{\stackrel{F_{k}\in \F}{F_{k-1}\subset F_{k}}}
\frac{1}{{n-|F_{k-1}| \choose n-|F_{k}|}}\\
&& \label{eq:20}
+ \sum_{(F_{2}, \ldots, F_{k-1})\in \B}
\!\!\!\!\!\!\!
\frac{|F_{2}|! \cdots (n-|F_{k-1}|)!}{n!}
\sum_{\stackrel{F_1\in \F}{F_1\subset F_2}} \frac{1}{{|F_2| \choose |F_1|}}
\!\!\!\!
\sum_{\stackrel{F_{k}\in \F}{F_{k-1}\subset F_{k}}}
\frac{1}{{n-|F_{k-1}| \choose n-|F_{k}|}}.
\end{eqnarray}

For the summation over $\A$,
the number of $F_1$ satisfying $F_1\subset F_2$ is at most $s-1$. We have
\begin{equation}
  \label{eq:21}
\sum_{\stackrel{F_1\in \F}{F_1\subset F_2}} \frac{1}{{|F_2| \choose |F_1|}}
\leq \frac{(s-1)}{\frac{n}{2}-2\sqrt{n\ln n}}.
\end{equation}
Apply inequality (\ref{eq:21}) to the first summation in
(\ref{eq:20}).
\begin{eqnarray}
\nonumber
\hspace*{2cm}&& \sum_{(F_{2}, \ldots, F_{k-1})\in \A}
\!\!\!\!\!\!\! \frac{|F_{2}|! \cdots (n-|F_{k-1}|)!}{n!}
\sum_{\stackrel{F_1\in \F}{F_1\subset F_2}} \frac{1}{{|F_2|
\choose |F_1|}} \!\!\!\! \sum_{\stackrel{F_{k}\in
\F}{F_{k-1}\subset F_{k}}}
\frac{1}{{n-|F_{k-1}| \choose n-|F_{k}|}} \\
\nonumber
&\leq& \sum_{(F_{2}, \ldots, F_{k-1})\in \A}
\!\!\!\!\!\!\!
\frac{|F_{2}|! \cdots (n-|F_{k-1}|)!}{n!}
\!\!\!\!
\sum_{\stackrel{F_{k}\in \F}{F_{k-1}\subset F_{k}}}
\frac{1}{{n-|F_{k-1}| \choose n-|F_{k}|}}
\frac{(s-1)}{\frac{n}{2}-2\sqrt{n\ln n}}\\
\nonumber
&\leq& \sum_{\stackrel{F_{2}, \ldots, F_{k-1}\in \F}
{F_{2}\subset \cdots \subset F_{k-1}}}
\!\!\!\!\!\!\!
\frac{|F_{2}|! \cdots (n-|F_{k-1}|)!}{n!}
\!\!\!\!
\sum_{\stackrel{F_{k}\in \F}{F_{k-1}\subset F_{k}}}
\frac{1}{{n-|F_{k-1}| \choose n-|F_{k}|}}
\frac{(s-1)}{\frac{n}{2}-2\sqrt{n\ln n}}\\
\label{eq:22}
&=& \E{X\choose k-1} \frac{(s-1)}{\frac{n}{2}-2\sqrt{n\ln n}}.
\end{eqnarray}

For the summation over $\B$,
the number of $F_{k}$ satisfying $F_{k-1}\subset F_{k}$
is at most $t-1$.
We have
\begin{equation}
  \label{eq:23}
\sum_{\stackrel{F_k\in \F}{F_{k-1}\subset F_k}}
\frac{1}{{n-|F_{k-1}| \choose n-|F_{k}|}} \leq
\frac{(t-1)}{\frac{n}{2}-2\sqrt{n\ln n}}.
\end{equation}
An inequality similar to (\ref{eq:22}) can be obtained:
\begin{equation}
  \label{eq:24}
\sum_{(F_{2}, \ldots, F_{k-1})\in \B}
\!\!\!\!\!\!\!
\frac{|F_{2}|! \cdots (n-|F_{k-1}|)!}{n!}
\sum_{\stackrel{F_1\in \F}{F_1\subset F_2}} \frac{1}{{|F_2| \choose |F_1|}}
\!\!\!\!
\sum_{\stackrel{F_{k}\in \F}{F_{k-1}\subset F_{k}}}
\frac{1}{{n-|F_{k-1}| \choose n-|F_{k}|}}\leq
 \E{X\choose k-1}\frac{(t-1)}{\frac{n}{2}-2\sqrt{n\ln n}}.
\end{equation}

Combining inequalities (\ref{eq:20}), (\ref{eq:22}) and
(\ref{eq:24}), we have
\begin{equation}
  \label{eq:25}
\E{X\choose k} \leq \E{X\choose k-1} \frac{s+t-2}{\frac{n}{2}-2\sqrt{n\ln n}}.
\end{equation}
>From inequalities (\ref{eq:17}) and (\ref{eq:25}), and the fact
that $\E{X\choose k-1}>0$,
 we have
\[\frac{\epsilon}{k} < \frac{s+t-2}{\frac{n}{2}-2\sqrt{n\ln n}}\],
which contradicts our choice of $\epsilon$.
\hfill $\square$


{\bf Proof of Theorem \ref{t3}}:
Let $\F$ be a $T$-free family of subsets of $[n]$.
By removing at most $\frac{2^{n+1}}{n^2}$ subsets,
without loss of generality,
we can assume $\F$ consists of subsets
of sizes in
$(\frac{n}{2}- 2\sqrt{n\ln n},
\frac{n}{2}+ 2\sqrt{n\ln n})$
 and $|\F|> (1+\epsilon){n\choose \lfloor \frac{n}{2}\rfloor}$.
Here $\epsilon = \frac{2t}{n} + \frac{16t}{n\sqrt{n\ln n}}$.

Let $X$ be the same variable as defined in the proof of Theorem \ref{t1}.
Recall
\begin{equation}
  \label{eq:26}
  \E(X)= \sum_{F\in \F}\frac{1}{{n\choose |F|}}.
\end{equation}
We have
\begin{eqnarray}
\nonumber
\E(X)&=& \sum_{F\in \F}\frac{1}{{n\choose |F|}}\\
\nonumber
 &\geq& \frac{|\F|}{{n\choose \lfloor \frac{n}{2}\rfloor}}\\
\label{eq:27}
 &>& 1+\epsilon.
\end{eqnarray}
Using that the variance of $X$ is nonnegative (or applying
Lemma~\ref{l3} with $r=1$ and $k=2$) we have
\begin{equation}
\label{eq:28}
  \E{X\choose 2}\geq \frac{1}{2}\E(X)(\E(X)-1).
\end{equation}
>From inequality (\ref{eq:27}) and (\ref{eq:28}), we get
\begin{equation}
\label{eq:29}
 \E{X\choose 2}
> \frac{\epsilon}{2}\E(X).
\end{equation}
A simple case of inequality (\ref{eq:18}) with $k=2$ is
\begin{equation}
  \label{eq:30}
\E{X\choose 2}= \sum_{\stackrel{F_1,F_2\in \F}{F_1\subset F_2}}
\frac{|F_1|!(|F_2|-|F_1|)!(n-|F_2|)!}{n!}.
\end{equation}

Now partition $\F$ into $\A\cup \B$ randomly. With probability
$\frac{1}{4}$, a pair $(F_1,F_2)$ has $F_1\in \A$ and $F_2\in \B$.
There is a partition $\F=\A\cup \B$ satisfying
\begin{equation}
  \label{eq:31}
\sum_{\stackrel{F_1\in \A,F_2\in \B}{F_1\subset F_2}}
\frac{|F_1|!(|F_2|-|F_1|)!(n-|F_2|)!}{n!}
>  \frac{\epsilon}{8}\E(X).
\end{equation}
Now we consider an edge-weighted
 bipartite graph $G$ with $V(G)=\A\cup \B$.
such that $F_1F_2$ is an edge of $G$ if $F_1\in \A$, $F_2\in \B$,
and $F_1\subset F_2$. Each edge $F_1F_2$ has weight
$\frac{|F_1|!(|F_2|-|F_1|)!(n-|F_2|)!}{n!}$. Inequality
(\ref{eq:31}) states that the total sum of edge-weights is greater
than $\frac{\epsilon}{8}\E(X)$.

 For any $F_1\in \A$, the weighted degree of $F_1$  is
\begin{equation}
  \label{eq:32}
d_{F_1}= \frac{1}{{n\choose |F_1|}}
\sum_{\stackrel{F_2\in \B}{F_1\subset F_2}}
\frac{1}{{n-|F_1|\choose n-|F_2|}}.
\end{equation}
Similarly, the weighted degree of $F_2\in \B$ is
\begin{equation}
  \label{eq:33}
d_{F_2}=\frac{1}{{n\choose |F_2|}}\sum_{\stackrel{F_1\in \A}{F_1\subset F_2}}
\frac{1}{{|F_2|\choose |F_1|}}.
\end{equation}

We delete vertices $F$ with weighted degree less than
$\frac{\epsilon}{8}\frac{1}{{n\choose |F|}}$ recursively until all
remaining vertices have weighted degree at least
$\frac{\epsilon}{8}\frac{1}{{n\choose |F|}}$ in the remaining
graph, call it $G'$, which has vertex partition $\A'\cup \B'$ with
$A'\subset A$ and $B'\subset B$. The sum of edge-weights in $G'$
is at least
\begin{eqnarray*}
\hspace*{2cm}&& \sum_{\stackrel{F_1\in \A',F_2\in \B'}{F_1\subset
F_2}}
\frac{|F_1|!(|F_2|-|F_1|)!(n-|F_2|)!}{n!} \\
&\geq& \sum_{\stackrel{F_1\in \A,F_2\in \B}{F_1\subset F_2}}
\frac{|F_1|!(|F_2|-|F_1|)!(n-|F_2|)!}{n!} \\
&&- \sum_{\stackrel{F_1\in \A\setminus \A',F_2\in \B}{F_1\subset F_2}}
\frac{|F_1|!(|F_2|-|F_1|)!(n-|F_2|)!}{n!} \\
&&-\sum_{\stackrel{F_1\in \A,F_2\in \B\setminus \B'}{F_1\subset F_2}}
\frac{|F_1|!(|F_2|-|F_1|)!(n-|F_2|)!}{n!}\\
&=&\sum_{\stackrel{F_1\in \A,F_2\in \B}{F_1\subset F_2}}
\frac{|F_1|!(|F_2|-|F_1|)!(n-|F_2|)!}{n!}\\
&&-\sum_{F_1\in \A\setminus \A'}\frac{d_{F_1}}{{n\choose |F_1|}}
-\sum_{F_2\in \B\setminus \B'}\frac{d_{F_2}}{{n\choose |F_2|}}\\
&>& \sum_{\stackrel{F_1\in \A,F_2\in \B}{F_1\subset F_2}}
\frac{|F_1|!(|F_2|-|F_1|)!(n-|F_2|)!}{n!}\\
&&- \sum_{F_1\in \A\setminus \A'}\frac{\epsilon}{8{n\choose|F_1|}}
-\sum_{F_2\in \B\setminus \B'}\frac{\epsilon}{8{n\choose |F_2|}}\\
&\geq& \sum_{\stackrel{F_1\in \A,F_2\in \B}{F_1\subset F_2}}
\frac{|F_1|!(|F_2|-|F_1|)!(n-|F_2|)!}{n!}  - \sum_{F\in \F}
\frac{\epsilon}{8{n\choose |F|}}\\
&\geq&\sum_{\stackrel{F_1\in \A,F_2\in \B}{F_1\subset F_2}}
\frac{|F_1|!(|F_2|-|F_1|)!(n-|F_2|)!}{n!}-\frac{\epsilon}{8}\E(X).
\end{eqnarray*}
Since the last expression is positive by (\ref{eq:31}), both
families $\A'$ and $\B'$ are non-empty.

By construction, every vertex in the remaining bipartite graph
$G'$ has weighted degree at least
$\frac{\epsilon}{8}\frac{1}{{n\choose |F|}}$. For any $F_1\in
\A'$, by (\ref{eq:32}) we have
\begin{equation}
  \label{eq:34}
\sum_{\stackrel{F_2\in \B}{F_1\subset F_2}}
\frac{1}{{n-|F_1|\choose |F_2|-|F_1|}}\geq \frac{\epsilon}{8}.
\end{equation}
Note
\begin{equation}
  \label{eq:35}
{n-|F_1|\choose |F_2|-|F_1|}\geq n-|F_1|
\geq \frac{n}{2}-2\sqrt{n\ln n}.
\end{equation}
Combining inequalities (\ref{eq:34}) and (\ref{eq:35}), we have
\begin{equation}
  \label{eq:36}
\sum_{\stackrel{F_2\in \B'}{F_1\subset F_2}} 1\geq
\frac{\epsilon}{8}(\frac{n}{2}-2\sqrt{n\ln n}).
\end{equation}
Similarly, for any $F_1\in \A'$,
\begin{equation}
  \label{eq:37}
\sum_{\stackrel{F_1\in \A'}{F_1\subset F_2}} 1\geq
\frac{\epsilon}{8}\left(\frac{n}{2}-2\sqrt{n\ln n}\right).
\end{equation}
In other words, the minimum degree (in the usual sense)
of $G'$ is at least $\frac{\epsilon}{8}(\frac{n}{2}-2\sqrt{n\ln n})>t$
for the choice of $\epsilon$.

A subgraph of $G'$ which is isomorphic to $T$ can be constructed
as follows. For any $u\in V(T)$, map $u$ to any vertex $v$ of $G'$.
Map the neighbors of $u$ in $T$ to the neighbors of $v$ in $G'$, and
so on. Since the minimum degree is at least $t$, we can always find
new vertex which has not been selected yet. This greedy algorithm
finds a subposet isomorphic to $T$.
\hfill $\square$

{\bf Proof of Theorem \ref{t5}:} Let $\F$ be any $P(G)$-free
subsets of $[n]$. By removing at most $\frac{2^{n+1}}{n^2}$
subsets, we can assume that $\F$ contains the subsets of sizes
only in the interval $(\frac{n}{2}-2\sqrt{n\ln
n},\frac{n}{2}+2\sqrt{n\ln n})$. Let $X$ be the random number
defined in the proof of Theorem \ref{t1}. We claim
$\E(X)=1+o_n(1)$.  Recall
\begin{equation}
  \label{eq:38}
\E(X)=\sum_{F\in \F}\frac{1}{{n\choose |F|}},
\end{equation}
so that $|\F|\le E(X)\nchn$.  We obtain an upper bound on $E(X)$.
As before,  we have
\begin{equation}
  \label{eq:39}
\E{X\choose 2} \geq \frac{1}{2}\E(X)(\E(X)-1).
\end{equation}

We will bound $\E{X\choose 2}$ in terms of $E(X)$. Recall
\begin{equation}
  \label{eq:40}
\E{X\choose 2}=\sum_{\stackrel{A,B\in \F}{A\subset B}}
\frac{|A|! (|B|-|A|)!(n-|B|)!}{n!}.
\end{equation}
We split the summation into two parts, depending on
whether $|B|-|A|=1$ or $|B|-|A|>1$.

For the case that $|B|-|A|>1$, let $Y$ be the random variable
counting a triple $(A,S,B)$ satisfying
\[A\subset S\subset B\quad A,B\in \F.\]
We have
\begin{eqnarray}
\nonumber
  \E(Y)&=& \sum_{\stackrel{A,B\in \F, S}{A\subset S\subset B}}
\frac{|A|! (|S|-|A|)!(|B|-|S|)!(n-|B|)!}{n!} \\
\nonumber
&=& \sum_{\stackrel{A,B\in \F}{A\subset B}}
\frac{|A|! (|B|-|A|)!(n-|B|)!}{n!} \sum_{S\colon A\subset S\subset B}
 \frac{1}{{|B|-|A|\choose |S|-|A|} }\\
\nonumber
&=& \sum_{\stackrel{A,B\in \F}{A\subset B}}
\frac{|A|! (|B|-|A|)!(n-|B|)!}{n!} (|B|-|A|-1) \\
\label{eq:41}
&\geq&  \sum_{\stackrel{A,B\in \F}{A\subset B, |B|-|A|>1}}
\frac{|A|! (|B|-|A|)!(n-|B|)!}{n!}.
\end{eqnarray}

Denote the number of vertices in $G$ by $v$  and the number of
edges in $G$ by $m$.  Since $\F$ is $P(G)$-free, there are no
$v+m$ subsets $A_1,A_2,\ldots, A_v, B_1,\ldots, B_m\in \F$
satisfying $A_i\subset S\subset B_j$ for $1\leq i\leq v$ and
$1\leq j\leq m$.

For any fixed subset $S$, either ``at most $m-1$ subsets in $\F$
are supersets of $S$'' or ``at most $v-1$ subsets in $\F$ are
subsets of $S$''. Define
\[\G_1=\{ S\mid |S|\in
(\frac{n}{2}-2\sqrt{n\ln n},\frac{n}{2}+2\sqrt{n\ln n}),
\mbox{$S$  has at most
$v-1$ subsets in $\F$}\}.\]
\[\G_2=\{ S\mid |S|\in
(\frac{n}{2}-2\sqrt{n\ln n},\frac{n}{2}+2\sqrt{n\ln n}),
\mbox{$S$  has at most
$m-1$ supersets in $\F$}\}.\]
$\G_1\cup \G_2$ covers all subsets with sizes in
$(\frac{n}{2}-2\sqrt{n\ln n},\frac{n}{2}+2\sqrt{n\ln n})$.
Rewrite $\E(Y)$ as
\begin{equation}
  \label{eq:42}
\E(Y) =
\sum_{S\colon ||S|-\frac{n}{2}|<2\sqrt{n\ln n}}
\frac{1}{{n\choose |S|}}
 \sum_{\stackrel{A\in \F}{A\subset S}}
\frac{1}{{|S|\choose |A|}}
 \sum_{\stackrel{B\in \F}{S\subset B}}
\frac{1}{{n-|S|\choose n-|B|}}.
\end{equation}

For $S\in\G_1$, we have
\begin{equation}
  \label{eq:43}
\sum_{B\in \F, S\subset B}\frac{1}{{n-|S|\choose n-|B|}}
\leq \frac{m-1}{\frac{n}{2}-2\sqrt{n\ln n}}.
\end{equation}
It implies
\begin{eqnarray}
\nonumber
  \sum_{S\in \G_2}
\frac{1}{{n\choose |S|}}
 \sum_{\stackrel{A\in \F}{A\subset S}}
\frac{1}{{|S|\choose |A|}}
 \sum_{\stackrel{B\in \F}{S\subset B}}
\frac{1}{{n-|S|\choose n-|B|}}
&\leq&
  \sum_{S\in \G_1}
\frac{1}{{n\choose |S|}}
 \sum_{\stackrel{A\in \F}{A\subset S}}
\frac{1}{{|S|\choose |A|}}
 \frac{m-1}{\frac{n}{2}-2\sqrt{n\ln n}}\\
\label{eq:44}
&\leq&
\E(X) 4\sqrt{n\ln n}\frac{m-1}{\frac{n}{2}-2\sqrt{n\ln n}}.
\end{eqnarray}

Similarly, we have
\begin{equation}
  \label{eq:45}
  \sum_{S\in \G_2}
\frac{1}{{n\choose |S|}}
 \sum_{\stackrel{A\in \F}{A\subset S}}
\frac{1}{{|S|\choose |A|}}
 \sum_{\stackrel{B\in \F}{S\subset B}}
\frac{1}{{n-|S|\choose n-|B|}}
\leq
\E(X) 4\sqrt{n\ln n}\frac{v-1}{\frac{n}{2}-2\sqrt{n\ln n}}.
\end{equation}
Combining equality (\ref{eq:42}) with inequalities  (\ref{eq:44})
and (\ref{eq:45}), we have
\begin{equation}
\E(Y)\leq \E(X) 4\sqrt{n\ln n}\frac{v+m-2}{\frac{n}{2}-2\sqrt{n\ln n}}.
\end{equation}

In particular, combining with inequality (\ref{eq:41}), we have
\begin{equation}
  \label{eq:46}
\sum_{\stackrel{A,B\in \F, |B|-|A|>1}{A\subset B}}
\frac{|A|! (|B|-|A|)!(n-|B|)!}{n!}
\leq \E(X) 4\sqrt{n\ln n}\frac{v+m-2}{\frac{n}{2}-2\sqrt{n\ln n}}
=o_n(\E(X)).
\end{equation}

Now we consider pairs $(A,B)$ with additional property $|B|-|A|=1$.
For any subset $S$, we define
\begin{eqnarray*}
N^+(S)&=&\{T\in \F\mid  S\subset T, |T|=|S|+1\}\\
N^-(S)&=&\{T\in \F\mid  T\subset S, |T|=|S|-1\}.
\end{eqnarray*}
Let $d^+(S)=|N^+(S)|$ and $d^-(S)=|N^-(S)|$.
We have
\begin{eqnarray}
\nonumber
\sum_{\stackrel{A,B\in \F}{A\subset B,|B|-|A|=1}}
\frac{|A|! (|B|-|A|)!(n-|B|)!}{n!}
&=&\sum_{\stackrel{A,B\in \F}{A\subset B,|B|-|A|=1}}
\frac{|A|! (n-|B|)!}{n!}\\
\label{eq:47}
&=&\sum_{A\in \F}\frac{d^+(A)}{{n\choose |A|}(n-|A|)}\\
\label{eq:48}
&=&\sum_{B\in \F}\frac{d^-(B)}{{n\choose |B|}|B|}.
\end{eqnarray}
We will show most contributions to the summation above are from
pairs $(A,B)$ with $d^+(A)\geq m$ and $d^-(A)\leq v$. We define
two subfamilies of $\F$ as follows:
\begin{eqnarray*}
\F_1&=&\{S\in \F\mid d^+(S)\geq  m\}\\
\F_2 &=&\{S\in \F\mid d^-(S)\geq  v\}.
\end{eqnarray*}

We have
\begin{eqnarray}
\nonumber
\hspace*{2cm}&&\hspace*{-3cm}
 \sum_{\stackrel{A,B\in \F}{A\subset B,|B|-|A|=1}}
\frac{|A|! (n-|B|)!}{n!} \\
\nonumber
&\leq&  \sum_{\stackrel{A\in \F_1, B\in \F_2}{A\subset B, |B|-|A|=1}}
\frac{|A|! (n-|B|)!}{n!}
+ \sum_{A\in \F\setminus \F_1}\frac{d^+(A)}{{n\choose |A|}(n-|A|)}
+\sum_{B\in \F\setminus \F_2}\frac{d^-(B)}{{n\choose |B|}|B|}\\
\nonumber
&\leq&\sum_{\stackrel{A\in \F_1, B\in \F_2}{A\subset B, |B|-|A|=1}}
\frac{|A|! (n-|B|)!}{n!}
+ \sum_{A\in \F\setminus \F_1}\frac{m-1}{{n\choose |A|}(n-\sqrt{2n\ln n})}\\
\nonumber
&&+\sum_{B\in \F\setminus \F_2}\frac{v-1}{{n\choose |B|}(n-\sqrt{2n\ln n})}\\
\label{eq:49}
&\leq&\sum_{\stackrel{A\in \F_1, B\in \F_2}{A\subset B, |B|-|A|=1}}
\frac{|A|! (n-|B|)!}{n!}
+ \frac{v+m-2}{n-\sqrt{2n\ln n}} \E(X).
\end{eqnarray}

Recall $C_{\sigma}$ is a random full chain of subsets of $[n]$.
For $i=1,2$, let $X_i=|\F_i\cap C_\sigma|$, so that
\begin{equation}
  \label{eq:50}
\E(X_i)=\sum_{F\in \F_i}\frac{1}{{n \choose |F|}}.
\end{equation}

Since $\F$ is $P(G)$-free, we have
$\F_1\cap \F_2=\emptyset$. In particular,
\begin{equation}
  \label{eq:51}
\E(X_1)+\E(X_2)\leq \E(X).
\end{equation}

Let us consider a ``diamond'' configuration $S\subset A_i\subset
B$ for ($i=1,2$) with $A_1, A_2\in \F_1$, $B\in \F_2$, and
$|B|=|S|+2$. In other words,
 $S=A_1\cap A_2$ and $B=A_1\cup A_2\in F_2$
where $A_1$ and $A_2$ ($\in\F_1$) only differ by one element. For
a fixed $S$, we define an auxiliary graph $L_S$ with vertex set
$N^+(S)\cap \F_1$ such that two subsets $A_1,A_2$ form an edge in
$L_S$ if $A_1\cup A_2\in \F_2$. We have
\begin{enumerate}
\item $L_S$ is $G$-free since $\F$ is $P(G)$-free.
\item Each edge of $L_S$ is in one-to-one correspondence with
a diamond configuration as above.
\end{enumerate}

Recall that the Tur\'an number $t(n,G)$ is the maximum number of
edges that a graph on $n$ vertices can have without containing the
subgraph $G$. The Erd\H{o}s-Simonovits-Stone Theorem
\cite{erdossimonovits, erdosstone} states
\begin{equation}
  \label{eq:52}
t(n,G)=\left(1-\frac{1}{\chi(G)-1}+o_n(1)\right)\frac{n^2}{2}.
\end{equation}
where $\chi(G)$ is the chromatic number of $G$.

Let $d^+_1(S)=|N^+(S)\cap \F_1|$ and $d^-_2(B)=|N^-(B)\cap \F_2|$.
The number of edges in $L_S$ is at most $t(d^+_1(S),G)$. We have
\begin{equation}
  \label{eq:53}
\sum_S f(|S|) t(d^+_1(S),G)\geq
\sum_{B\in \F}f(|B|-2) {d^-_2(B)\choose 2}.
\end{equation}
Here $f(k)$ is any nonnegative function over integers and the
summation on the left is taken over all $S$ with sizes in
$(\frac{n}{2}-2\sqrt{n\ln n}-1,\frac{n}{2}+2\sqrt{n\ln n}-1)$.
Choose $f(k)=\frac{1}{{n\choose k} (n-k)^2}$ for $k\in
(\frac{n}{2}-2\sqrt{n\ln n}-1,\frac{n}{2}+2\sqrt{n\ln n}-1)$. We
have
\begin{eqnarray}
\nonumber
  \sum_S f(|S|) t(d^+_1(S),G) &=&
\left(1-\frac{1}{\chi(G)-1}+o_n(1)\right)
\sum_S f(|S|) \frac{(d^+_1(S))^2}{2}\\
\nonumber &\leq&
\frac{1}{2}\left(1-\frac{1}{\chi(G)-1}+o_n(1)\right)
\sum_S f(|S|) d^+_1(S)(n-|S|)\\
\nonumber &=& \frac{1}{2}\left(1-\frac{1}{\chi(G)-1}+o_n(1)\right)
\sum_S \frac{d^+_1(S)}{{n\choose |S|}(n-|S|)}\\
\label{eq:54} &=&
\frac{1}{2}\left(1-\frac{1}{\chi(G)-1}+o_n(1)\right)\E(X_1).
\end{eqnarray}
\begin{eqnarray}
\nonumber
   \sum_{B\in \F_2}f(|B|-2) {d^-_1(B)\choose 2}
&=&\frac{1}{2} \sum_{B\in \F_2} \frac{1}{{n\choose |B|-2}(n-|B|+2)^2}
{(d^-_2(B))^2-d^-_2(B)} \\
\nonumber &=&\frac{1}{2}\left(1+O(\frac{\sqrt{n\ln n}}{n})\right)
\sum_{B\in \F_2} \frac{(d^-_2(B))^2-d^-_2(B)}{{n\choose |B|}|B|^2}\\
\label{eq:55} &=&\frac{1}{2}\left(1+O(\frac{\sqrt{n\ln
n}}{n})\right) \sum_{B\in \F_2} \frac{(d^-_2(B))^2}{{n\choose
|B|}|B|^2} -O(\frac{1}{n})\E(X_2).
\end{eqnarray}
Applying the Cauchy-Schwartz Inequality. the inequalities above,
and the Arithmetic-Geometric Mean Inequality, we have
\begin{eqnarray}
\nonumber
\sum_{\stackrel{A\in \F_1, B\in \F_2}{A\subset B, |B|-|A|=1}}
\frac{|A|! (n-|B|)!}{n!}
&=&
\sum_{B\in \F_2} \frac{d^-_2(B)}{{n\choose |B|}|B|} \\
\nonumber &\leq& \sqrt{{\sum_{B\in \F_2} \frac{1}{{n\choose |B|}}}
 {\sum_{B\in \F_2} {\frac{(d^-_B)^2}{{n\choose |B|}}|B|^2}}}\\
\nonumber
&\leq & \sqrt{\E(X_2)
\left(1-\frac{1}{\chi(G)-1}+o_n(1)\right)\E(X_1)}\\
\nonumber
&=& \left(\sqrt{1-\frac{1}{\chi(G)-1}} +o_n(1)\right)
\sqrt{\E(X_1)\E(X_2)}\\
\nonumber
&\leq& \left(\sqrt{1-\frac{1}{\chi(G)-1}} +o_n(1)\right)
\frac{\E(X_1)+\E(X_2)}{2}\\
\label{eq:56} &\leq& \left(\sqrt{1-\frac{1}{\chi(G)-1}}
+o_n(1)\right)\frac{\E(X)}{2}.
\end{eqnarray}
Combining inequalities (\ref{eq:46}), (\ref{eq:49}), and
(\ref{eq:56}) , we have
\begin{eqnarray}
\nonumber
  \E{X\choose 2} &=& \sum_{\stackrel{A,B\in \F, |B|-|A|>1}{A\subset B}}
\frac{|A|! (|B|-|A|)!(n-|B|)!}{n!}\\
\nonumber
&&+ \sum_{\stackrel{A,B\in \F, |B|-|A|=1}{A\subset B}}
\frac{|A|! (|B|-|A|)!(n-|B|)!}{n!}\\
\nonumber &\leq& o_n(\E(X)) + \left(\sqrt{1-\frac{1}{\chi(G)-1}}
+o_n(1)\right)
\frac{\E(X)}{2}\\
\label{eq:57} &\leq& \left(\sqrt{1-\frac{1}{\chi(G)-1}}
+o_n(1)\right)\frac{1}{2}\E(X).
\end{eqnarray}
Combining inequalities (\ref{eq:39}) and (\ref{eq:57}), we have
\begin{equation}
  \label{eq:58}
E(X)\leq 1 + \sqrt{1-\frac{1}{\chi(G)-1}} +o_n(1).
\end{equation}
The proof is finished by observing $|\F|\leq \E(X)
{n\choose \lfloor \frac{n}{2}\rfloor}.$ \hfill $\square$

\section{Further research}

Let $$\pi(H):=\lim_{n\rightarrow\infty} \frac{\La(n,H)}{\nchn},$$
when this limit exists.  Does this limit exist for all posets $H$,
and, if so, how does it depend on $H$?  For posets $H$ of height
two, we know that the limit, when it exists, belongs to the
interval $[1,2]$. Are there any $H$ of height two such that
$\pi(H)$ is strictly between 1 and 2?

More generally, for all posets $H$ where we know $\pi(H)$,
$\pi(H)$ is an integer.  Is this true in general?  In fact,
examples we looked have have $\pi(H)$ equal to the maximum number
$m$ such that the middle $m$ levels of the Boolean lattice
$B_n=(2^{[n]},\subseteq)$ do not contain $H$, no matter how large
$n$ is (as observed by Mike Saks and Pete Winkler, unpublished).

We once asked whether there exists a number $c_h$ such that for
all posets $H$ of height $h$, $\pi(H)\le c_h$. As we noted above,
$c_2=2$. However, Lu and, independently, Tao Jiang, pointed out
that no such $c_h$ for $h\ge3$.  The idea is that if one takes
$\F$ to consist of the middle $m$ levels in the Boolean lattice
$B_n$, then two sets $A,B\in \F$ with $A\subset B$ have at most
$2^{m-1}-2$ sets $C$ with $A\subset C\subset B$.  Hence, the
family $\F$, which has size $\sim m\nchn$, avoids the height 3
poset consisting of a minimum element, a maximum element, and an
antichain of $2^{m-1}-1$ elements in between.  This forces $c_3$
to be larger than any $m$, so that no such $c_3$ exits.  It seems
that not just the height, but the width, of $H$ affects $\pi(H)$.

It would therefore be interesting to determine $\pi(B_n)$ for the
Boolean lattice $B_n$ . The smallest crown for which $\pi$ is not
yet determined is $\Oh_6$, the height two poset formed by the
middle two levels of $B_3$.  Even for a poset as fundamental as
the diamond poset $B_2$, we only know that $\pi(B_2)$, if it
exists, must be in the interval $[2,3]$.


\comments{
\section{Appendix}
\begin{lemma}
Suppose $X$ is a random variable which takes on nonnegative
integer values. Let $f(x)$ and $g(x)$ be two increasing functions.
Then}

\begin{thebibliography}{99}
\bibitem{bukh}
B. Bukh, {Set families with a forbidden poset}, preprint (2008).

\bibitem{chernoff}
H. Chernoff, A note on an inequality involving the normal
distribution, {\it Ann. Probab.} {\bf 9} (1981), 533-535.

\bibitem{dk} A. De Bonis, G. O.H. Katona, {Largest families
without an $r$-fork}, {\it Order} {\bf 24} (2007), 181--191.

\bibitem{dks}
A. De Bonis, G. O.H. Katona, K. J. Swanepoel, {Largest family
without $A\cup B \subset C\cap D$}, {\it J. Combin. Theory (Ser.
A)} {\bf 111} (2005), 331-336.

\bibitem{erdos} P. Erd\H{o}s,
On a lemma of Littlewood and Offord,
{\it Bull. Amer. Math. Soc.} {\bf 51} (1945), 898-902.

\bibitem{erdossimonovits}
P. Erd\H{o}s and M. Simonovits,
A limit theorem in graph theory,
{\it Studia Sci. Math. Hungar.} {\bf 1} (1966), 51-57.

\bibitem{erdosstone}
P. Erd\H{o}s and A. H. Stone,
On the structure of linear graphs,
{\it Bull. Amer. Math. Soc.} {\bf 52} (1946), 1087-1091.

\bibitem{fkg}
C. M. Fortuin, P. N. Kasteleyn, and J. Ginibre, Correlation
inequalities for some partially ordered sets, {\it Comm. Math.
Physics}, {\bf 22} (1971), 89-103.

\bibitem{gk} J. R. Griggs, G. O. H. Katona, {No
four subsets forming an $N$}, {\it J. Combinatorial Theory (Ser.
A)} {\bf 115} (2008), 677--685.

\bibitem{sperner} E. Sperner,
{Ein Satz \"{u}ber Untermegen einer endlichen Menge},
{\it Math. Z.} {\bf 27} (1928), 544-548.

\bibitem{thanh} H. T. Thanh,
An extremal problem with excluded subposets in the Boolean lattice,
{\it Order} {\bf 15} (1998), 51-57.

\end{thebibliography}
\end{document}